\documentclass[reqno,10pt]{amsart}
\usepackage{graphicx}

\def\RR    {{R\!\!\!\!\!I~}}

\newcommand \Scal {\mathcal S}
\newcommand \Rcal {\mathcal R}
\newcommand \CC {\mathcal C} 
\def\ld    {\lambda}

\newcommand \HH  {\mathcal H}

\newcommand \WW  {\mathcal W}
 
\def\pt    {\partial}

\def\RR    {{R\!\!\!\!\!I~}}

\def\pt   {\partial}

\newcommand{\rf}[1]{(\ref{#1})}
\begin{document}
\title[The Riemann problem for the shallow water equations]
{The Riemann problem for the shallow water equations
with discontinuous topography}\footnotetext{To appear in: Communications in Mathematical Sciences.}

\author[P.G. L{\tiny e}Floch and M.D. Thanh]
{Philippe G.  L{\tiny e}Floch and  Mai Duc Thanh}

\address{ Philippe G. LeFloch\newline
Laboratoire Jacques-Louis Lions \& Centre National de la Recherche
Scientifique, Universit\'e de Paris 6, 4 Place Jussieu, 75252 Paris,
France.} \email{LeFloch@ann.jussieu.fr}

\address{Mai Duc Thanh \newline
Department of Mathematics, International University, Quarter 6,
Linh Trung Ward, Thu Duc District, Ho Chi Minh City, Vietnam}
\email{MDThanh@hcmiu.edu.vn}


\begin{abstract} We construct the solution of the Riemann problem for the shallow water equations with
discontinuous topography. The system under consideration is non-strictly hyperbolic 
and does not admit a fully conservative form, and we establish the existence of two-parameter
wave sets, rather than wave curves. The selection of admissible waves is particularly challenging. 
Our construction is fully explicit, and leads to formulas that can be implemented 
numerically for the approximation of the general initial-value problem. 
\end{abstract}

\maketitle \numberwithin{equation}{section}
\newtheorem{theorem}{Theorem}[section]
\newtheorem{lemma}[theorem]{Lemma}
\newtheorem{proposition}[theorem]{Proposition}
\newtheorem{definition}[theorem]{Definition}
\newtheorem{corollary}[theorem]{Corollary}

\section{Introduction}

\subsection{Shallow water equations} 

We consider the one-dimensional shallow water equations
\begin{equation}
\aligned &\pt_th + \pt_x(hu)  =  0,
\\
&\pt_t(hu) + \pt_x(h(u^2 + g\frac{h}{2}))  = - gh\pt_x a,
\\
&\pt_t a = 0,\\
\endaligned
\label{1.1}
\end{equation}
where $h$ denotes the height of
the water from the bottom to the surface, $u$ the velocity of the fluid,
$g$ the gravity constant, and $a$ the height of the river bottom from a given level. 
Following LeFloch \cite{LeFloch89} we supplement the first two balance laws for the fluid with the equation 
$\pt_ta=0$ corresponding to a fixed geometry. 
Adding the equation $\pt_t a=0$ allows us to view the shallow water equations (the first two equations in \rf{1.1}), 
which form a strictly hyperbolic system of balance laws in nonconservative form, 
as a non-strictly hyperbolic system of balance laws with a linearly degenerate characteristic field.

We are mainly interested in the case that $a$ is piecewise constant
$$
a(x)=\left\{\begin{array}{ll}a_L,\quad &x< 0,\\
a_R,\quad &x>0, \end{array}\right.
$$
where $a_L,a_R$ are two distinct constants.  The Riemann problem associated with \rf{1.1} is the
initial-value problem corresponding to the initial conditions of
\begin{equation}
(h,u,a)(x,0)= \left\{\begin{array}{ll}(h_L,u_L,a_L),\quad &x< 0,\\
(h_R,u_R,a_R),\quad &x>0. \end{array}\right. 
\label{1.2}
\end{equation}
Since $a$ is discontinuous, the system \rf{1.1} cannot be written
in a fully conservative form, and the standard notion of weak solutions for
hyperbolic systems of conservation laws does not apply. However, the
equations still make sense within the framework introduced in Dal Maso, LeFloch, and Murat
\cite{DalMasoLeFlochMurat}. (For a recent review see \cite{LeFloch02,LeFloch04}.)  

%------------------------------------------------------------------------------------------------------

\subsection{DLM generalized Rankine-Hugoniot relations}

Consider an elementary discontinuity propagating with the speed $\ld$ and satisfying the equations \rf{1.1}. 
Observe that the
Rankine-Hugoniot relation associated with the third equation in 
\rf{1.1} simply reads
\begin{equation}
-\ld [a]  =  0,
\label{1.3}
\end{equation}
where $[a] := a_+ - a_-$ denotes the jump of the bottom level function
$a$, and $a_\pm$ denotes its left- and right-hand traces. Then, we have the following possibilities: 
\begin{itemize}\item[(i)] either the component $a$ remains constant across the propagating discontinuity,
 \item[(ii)] or $a$ changes its levels across the  discontinuity and  the
 discontinuity is stationary, i.e., the speed $\ld$ vanishes.
\end{itemize}
This observation motivates us to define the admissible elementary waves of the system \rf{1.1}.
First of all, assume that the bottom level $a$ remains constant across a discontinuity; 
then, $a$ should be constant in a neighborhood of the discontinuity. Eliminating $a$ from \rf{1.1},
we obtain the following system of two conservation laws
\begin{equation}
\aligned &\pt_th + \pt_x(hu)  =  0,
\\
&\pt_t(hu) + \pt_x(h(u^2 + g\frac{h}{2}))  = 0,
\\
\endaligned
\label{1.4}
\end{equation}
Thus, the left- and right-hand states are related by the Rankine Hugoniot relations corresponding to \rf{1.4}
\begin{equation}
\aligned &-\ld[h] + [hu]  =  0,
\\
&-\ld[hu] + [h(u^2 + g\frac{h}{2})]  = 0,
\\
\endaligned
\label{1.5}
\end{equation}
where $[h] : = h_+ - h_+$, etc.

Second, suppose that the component $a$ is discontinuous so that the speed vanishes. Then, the solution is
independent of the time variable, and it is natural to search for a solution obtained as 
the limit of a sequence of time-independent smooth solutions of \rf{1.1}. (See below.)  

Suppose that $(x,t) \mapsto (h, u, a)$ is a smooth solution of
\rf{1.1}. Then, the system \rf{1.1} can be written in the following form, as a system of conservation laws for
the (now conservative) variables $(h,u,a)$:
\begin{equation}
\aligned &\pt_th + \pt_x(hu)  =  0,
\\
&\pt_tu + \pt_x\big({u^2\over 2} + g(h+a)\big)  = 0,
\\
&\pt_t a =0.\\
\endaligned
\label{1.6}
\end{equation}
Hence, time-independent solutions of \rf{1.1} satisfy
\begin{equation}
\aligned &(hu)'  =  0,
\\
&\big({u^2\over 2} + g(h+a)\big)'  = 0,
\\
\endaligned
\label{1.7}
\end{equation}
where the dash denotes the differentiation with respect to $x$.
Trajectories initiating from a given state $(h_0,u_0,a_0)$ are given by 
\begin{equation}
\aligned &hu  =  h_0u_0,
\\
&{u^2\over 2} + g(h+a)  = {u_0^2\over 2} + g(h_0+a_0).
\\
\endaligned
\label{1.8}
\end{equation}
It follows from \rf{1.8} that the trajectories of \rf{1.7} can
be expressed in the form $u=u(h)$, $a=a(h)$.  Now, letting $h\to
h_\pm$ and setting $u_\pm = u(h_\pm)$, $a_\pm = a(h_\pm)$, we see
that the states $(h_\pm, u_\pm, a_\pm)$ satisfy the
Rankine-Hugoniot relations associated with \rf{1.6}, but with zero shock speed:
\begin{equation}
\aligned &[hu]  =  0,
\\
&[{u^2\over 2} + g(h+a)]  = 0,
\\
\endaligned
\label{1.9}
\end{equation}

The above discussion leads us to define the elementary waves of interest, as follows.

\begin{definition}  The admissible waves for the system \rf{1.1} are the following ones:
\begin{itemize} \item[(a)] the {\rm rarefaction waves}, which
are smooth solutions of \rf{1.1} with constant component $a$
depending only on the self-similarity variable $x/t$;
\item[(b)]  the {\rm shock waves} which satisfy \rf{1.5} and Lax
shock inequalities and have constant component $a$; 

\item[(c)] and the {\rm stationary waves} which have zero speed and satisfy \rf{1.9}.
\end{itemize}\end{definition}

As will be checked later, the system \rf{1.1} is {\it not strictly hyperbolic}, as was already 
observed in the previous work \cite{LeFlochThanh03.2}. Recall that therein
we studied the Riemann problem in a nozzle with variable cross-section and constructed all of the Riemann
solutions. The present model is analogous, and our main purpose in the present paper is to 
demonstrate that the technique in  \cite{LeFlochThanh03.2} extends to the shallow water model 
and to construct the solution of the Riemann problem. The lack of strict hyperbolicity and the 
nonconservative form of the equation make the problem particularly challenging. 
Some aspects of this problem are also covered by Alcrudo and Benkhaldoun \cite{AB}. 
For works on various related models including scalar conservation laws 
we refer to \cite{MarchesinPaes-Leme,IsaacsonTemple95, IsaacsonTemple92,HouLeFloch,
HayesLeFloch,Gosse,GoatinLeFloch,AGG}. 

%---------------------------------------------------------------------------------------------------------------

\subsection{Results and perspectives}

As we will show, waves in the same characteristic field may be repeated in a single
Riemann solution. This happens when waves cross the boundary of the strictly hyperbolic regions  
and the order of characteristic speeds changes.
We will also show below that the Riemann problem may not always have a solution. 
The Riemann problem may admit exactly one, or two, or up to three distinct solutions for different ranges 
of left-hand and right-hand states. Thus, uniqueness does not hold for the Riemann problem, as
was already observed for the nozzle flow system. 

Each possible construction leads to a solution that
depends continuously on the left-hand and right-hand states. This is a direct consequence of 
the smoothness of the elementary wave curves; by the implicit function theorem,
the intermediate waves depend continuously on their left- or right-hand states as well 
as on the Riemann data. These results agree with \cite{LeFlochThanh03.2} 
which covered fluids in a nozzle with variable cross section. 

In the present model, the curve of stationary wave is strictly convex. 
To find stationary waves, one needs to determine the roots of a nonlinear equation 
(see the function $\varphi$ in \eqref{4.3}) which is convex and, therefore, can be
easily computed numerically. 
The Riemann solver derived in the present paper should be useful in combination with 
numerical methods for shallow water systems developed in 
 \cite{AndrianovWarnecke, AudusseBouchutBristeauKleinPerthame, GreenbergLerouxBarailleNoussair, GT, KroenerThanh04,CGGP} for which we refer to the lecture notes by Bouchut \cite{Bouchut}.

%*****************************************************************************************

\section{Background}

\subsection{Shallow water equations as a non-strictly hyperbolic system}

We now discuss the system \rf{1.1} in the nonconservative variables $U=(h,u,a)$.
From \rf{1.6} if follows that, for smooth solutions, \rf{1.1} is equivalent to 
\begin{equation}
\aligned &\pt_th + u\pt_xh + h\pt_xu  =  0,
\\
&\pt_tu + g\pt_xh + u\pt_x u + g\pt_x a = 0,
\\
&\pt_t a =0,\\
\endaligned
\label{2.1}
\end{equation}
which can be written in the nonconservative form
\begin{equation}
\pt_t U + A(U)\pt_x U=0,
\label{2.2}
\end{equation}
where the Jacobian matrix $A(U)$ is given by
$$
A(U)=\left(\begin{matrix} u&u&0\\
g&u&g\\
0&0&0\end{matrix}\right).
$$

The eigenvalues of $A$ are 
\begin{equation}
\ld_1(U):=u-\sqrt{gh}<\ld_2(U):=u+\sqrt{gh}, \quad \ld_3(U):=0,
\label{2.3}
\end{equation}
and corresponding eigenvectors can be chosen as
\begin{equation}
\aligned & r_1(U):=(h,-\sqrt{gh},0)^t,\quad
r_2(U):=(h,\sqrt{gh},0)^t,\\
&r_3(U):=(gh,-gu,u^2-gh)^t.\endaligned \label{2.4}
\end{equation}
We see that the first and the third characteristic fields may coincide:
$$
(\ld_1(U),r_1(U)) = (\lambda_3(U),r_3(U))
$$
on a hypersurface in the variables $(h,u,a)$, which can be identified as
\begin{equation}
\CC_+:=\{(h,u,a) | \quad u=\sqrt{gh}\}.
\label{2.5}
\end{equation}

Similarly, the second and the third characteristic fields may coincide:
$$
(\ld_2(U),r_2(U))=(\lambda_3(U),r_3(U))
$$
on a hypersurface in the variables $(h,u,a)$, which can be identified as
\begin{equation}
\CC_-:=\{(h,u,a) | \quad u=-\sqrt{gh}\}.
\label{2.6}
\end{equation}
The third eigenvalue $(\ld_3,r_3)$ is linearly degenerate, and we have
$$
-\nabla \ld_1(U)\cdot r_1(U) =\nabla\ld_2(U)\cdot r_2(U)={3\over
2}\sqrt{gh}\ne 0,\quad h> 0.
$$
Note also that the first and the second
characteristic fields $(\ld_1,r_1)$, $(\ld_2,r_2)$ are genuinely
nonlinear in the open half-space  $\{(h,u,a) |\quad h>0\}$.

It is convenient to set
$$
\CC=\CC_+\cup \CC_-=\{(h,u,a) |\quad  u^2-gh=0\},
$$
which is the hypersurface on which the system fails to be strictly hyperbolic.

\begin{figure} 
  \includegraphics[width=0.7\textwidth]{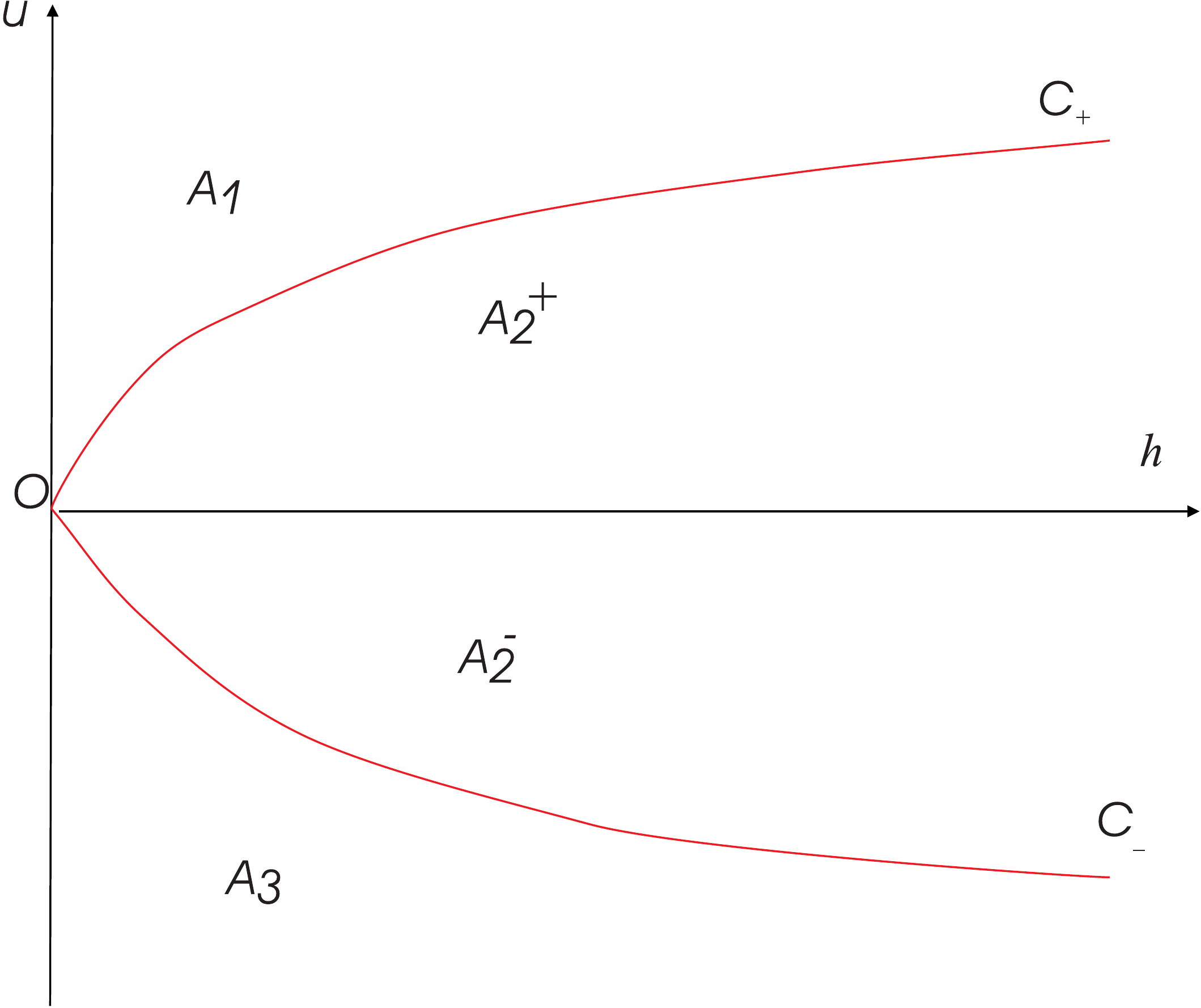}\\
  \caption{Projection of strictly hyperbolic regions in the $(h,u)$-plane}\label{fig21}
\end{figure}

In conclusion we have established (cf.~Figure \ref{fig21}):

\begin{lemma}
On the  hypersurface $\CC_+$ in the variables $(h,u,a)$ the first and the third characteristic speeds
coincide and, on the hypersurface $\CC_-$, the second and the third characteristic speeds coincide.
Hence, the system \rf{1.1} is non-strictly hyperbolic. 
\end{lemma}

The hypersurface $\CC$ divides the phase
domain into three disjoint regions, denoted by $A_1, A_2$ and $A_3$, in which the
system is strictly hyperbolic. More precisely, we define 
\begin{equation}
\aligned &A_1:=\{(h,u,a) \in \RR_+\times\RR\times\RR_+ |\quad \ld_2(U)>\ld_1(U)>\ld_3(U)\},\\
&A_2:=\{(h,u,a) \in \RR_+\times\RR\times\RR_+ | \quad \ld_2(U)>\ld_3(U)>\ld_1(U)\},\\
&A_2^+:=\{(h,u,a) \in A_2 | \quad u>0\},\\
&A_2^-:=\{(h,u,a) \in A_2 | \quad u<0\},\\
&A_3:=\{(h,u,a)\in \RR_+\times\RR\times\RR_+  |\quad \ld_3(U)>\ld_2(U)>\ld_1(U)\}.\\
\endaligned
\label{2.7}
\end{equation}
The strict hyperbolicity domain is not connected, which makes the Riemann
problem delicate to solve.

%==========================================================================================

\subsection{Wave curves}

We begin by investigating some properties of the curves of admissible waves. 

First, consider shock curves from a given left-hand state
$U_0=(h_0,u_0,a_0)$ consisting of all right-hand  states
$U=(h,u,a)$ that can be connected to $U_0$ by a shock wave. Thus,
it follows from \rf{1.6} that $U$ and $U_0$ are related by the
Rankine-Hugoniot relations
\begin{equation}
\aligned &-\bar\ld[h] + [hu]  =  0,
\\
&-\bar\ld[hu] + [h(u^2 + g\frac{h}{2}] = 0,
\\
\endaligned
\label{3.1}
\end{equation}
where $[h]=h-h_0$, etc, and $\bar\ld=\bar\ld(U_0,U)$ is the shock speed.

Fix the state $U_0$. 
A straightforward calculation from the Rankine-Hugoniot
relations \rf{3.1} shows that the restriction to the $(h,u)$ plane of 
the Hugoniot set consists of two curves given by 
\begin{equation}
u=u_0\pm \sqrt{g\over 2}(h-h_0)\sqrt{\Big({1\over h}+{1\over
h_0}\Big)}. \label{3.2}
\end{equation}
Moreover, along these two curves it holds
$$
\aligned 
{du\over dh} 
& = \pm \sqrt{g\over 2}\Big(\sqrt{{1\over h}+{1\over
h_0}} -(h-h_0){1\over 2h^2\sqrt{\dfrac{1}{h}+\dfrac{1}{h_0}}}\Big)
\\
& \to\pm\sqrt{g\over h_0} \text{ as $h\to h_0$.}
\endaligned 
$$
Since the $i$th-Hugoniot curve
is tangent to $r_i(U_0)$ at $U_0$, we conclude that the first
Hugoniot curve associated with the first characteristic field is
\begin{equation}
\HH_1(U_0) :\quad u:=u_1(h,U_0)=u_0-  \sqrt{g\over
2}(h-h_0)\sqrt{\Big({1\over h}+{1\over h_0}\Big)}, \quad h\ge
0,\label{3.3}
\end{equation}
while one associated with the second
characteristic field is 
\begin{equation}
\HH_2(U_0):\quad u:=u_2(h,U_0)=u_0+  \sqrt{g\over
2}(h-h_0)\sqrt{\Big({1\over h}+{1\over h_0}\Big)},\quad h\ge 0.
\label{3.4}
\end{equation}

Along the Hugoniot curves $\HH_1, \HH_2$, the corresponding shock
speeds are given by
\begin{equation}
\aligned
\bar\ld_{1,2}(U_0,U) &={hu_{1,2}-h_0u_0\over h-h_0}\\
&=u_0\mp\sqrt{{g\over 2}\Big(h+{h^2\over h_0}\Big)},\quad h\ge
0,\endaligned \label{3.5}
\end{equation}

As is customary, the shock speed $\bar\ld_i(U_0,U)$ is required to satisfy Lax
shock inequalities \cite{Lax71}: 
\begin{equation}
\ld_i(U)<\bar\ld_i(U_0,U)<\ld_i(U_0),\quad i=1,2. \label{3.6}
\end{equation}
Thus, the $1$-shock curve $\Scal_1(U_0)$ initiating
from the left-hand state $U_0$ and consisting of all right-hand states
$U$ that can be connected to $U_0$ by a Lax shock associated with
the first characteristic field is
\begin{equation}
\Scal_1(U_0) :\quad u=u_1(h,U_0)=u_0-  \sqrt{g\over
2}(h-h_0)\sqrt{\Big({1\over h}+{1\over h_0}\Big)}, \quad h>
h_0.
\label{3.7}
\end{equation}

Similarly, the $2$-shock curve $\Scal_2(U_0)$ issuing from a
left-hand state $U_0$ consisting of all right-hand states $U$ that
can be connected to $U_0$ by a Lax shock associated with the
second characteristic field is 
\begin{equation}
\Scal_2(U_0) :\quad u=u_1(h,U_0)=u_0+ \sqrt{g\over
2}(h-h_0)\sqrt{\Big({1\over h}+{1\over h_0}\Big)}, \quad h<
h_0,\label{3.8}
\end{equation}

We summarize these results in the following proposition.

\begin{lemma}[Shock wave curves] 
Given a left-hand state $U_0$, the $1$-shock curve $\Scal_1(U_0)$
consisting of all right-hand states $U$ that can be connected to
$U_0$ by a Lax shock is 
$$
\Scal_1(U_0) :\quad u=u_1(h,U_0)=u_0-  \sqrt{g\over
2}(h-h_0)\sqrt{\Big({1\over h}+{1\over h_0}\Big)}, \quad h> h_0. 
$$
The $2$-shock curve $\Scal_2(U_0)$ consisting of all right-hand
states $U$ that can be connected to $U_0$ by a Lax shock is 
$$
\Scal_2(U_0) :\quad u=u_1(h,U_0)=u_0+ \sqrt{g\over
2}(h-h_0)\sqrt{\Big({1\over h}+{1\over h_0}\Big)}, \quad h< h_0.
$$
\end{lemma}

In view of the Lax shock inequalities \rf{3.6}, we also conclude that 
the backward $1$-shock curve $\Scal_1^B(U_0)$
issuing from a right-hand state $U_0$ and consisting of all left-hand
states $U$ that can be connected to $U_0$ by a Lax shock
associated with the first characteristic field is 
\begin{equation}
\Scal_1^B(U_0) :\quad u=u_1(h,U_0)=u_0-  \sqrt{g\over
2}(h-h_0)\sqrt{\Big({1\over h}+{1\over h_0}\Big)}, \quad h<
h_0,\label{3.9}
\end{equation}
Similarly, the backward $2$-shock curve $\Scal_2^B(U_0)$ issuing
from  a right-hand state $U_0$ and consisting of all left-hand states
$U$ that can be connected to $U_0$ by a Lax shock associated with
the second characteristic field is 
\begin{equation}
\Scal_2^B(U_0) :\quad u=u_1(h,U_0)=u_0+ \sqrt{g\over
2}(h-h_0)\sqrt{\Big({1\over h}+{1\over h_0}\Big)}, \quad h>
h_0,\label{3.10}
\end{equation}

Next, we discuss the properties of rarefaction waves, i.e., smooth
self-similar solutions to the system \rf{1.1} associated with one of the two
genuinely nonlinear characteristic fields. These waves satisfy the
ordinary differential equation: 
\begin{equation}
\frac{dU}{d\xi} = \frac{r_i(U)}{\nabla \ld_i\cdot r_i(U)},\quad
\xi =x/t,\quad i=1,2.\label{3.11}
\end{equation}
For waves in the first family, we have
$$
\aligned
\frac{dh(\xi)}{d\xi} &= -\frac{2h(\xi)}{3\sqrt{gh(\xi)}}=-{2\over 3\sqrt{g}}\sqrt{h(\xi)},\\
\frac{du(\xi)}{d\xi} &= \frac{-2\sqrt{gh(\xi)}}{-3\sqrt{gh(\xi)}}={2\over 3},\\
\frac{da(\xi)}{d\xi} &= 0.\\
\endaligned
$$
It follows that
$$
\frac{du}{dh} = -\sqrt{\frac{g}{h}},
$$
therefore, the integral curve passing through a given point $U_0=(h_0,u_0,a_0)$ is given by
$$
u=u_0-2\sqrt{g}(\sqrt{h}-\sqrt{h_0}).
$$

Moreover, the characteristic speed should increase through a rarefaction fan, i.e.,
\begin{equation}
\ld_1(U)\ge \ld_1(U_0),\label{3.14}
\end{equation}
which implies
$$
h\ge h_0,
$$
Thus, we define a rarefaction curve $\Rcal_1(U_0)$
issuing from a given left-hand state $U_0$ and consisting of all the
right-hand states $U$ that can be connected to $U_0$ by a
rarefaction wave associated with the first characteristic field as
\begin{equation}
\Rcal_1(U_0):\quad
u=v_1(h,U_0):=u_0-2\sqrt{g}(\sqrt{h}-\sqrt{h_0}),\quad h\le
h_0.\label{3.16}
\end{equation}
A $1$-rarefaction wave is determined by
\begin{equation}
u=u_0+{2\over 3}\Big({x\over t}-{x_0\over t_0}\Big) \label{3.17}
\end{equation}
while $h$ is determined by the equation \rf{3.16} and the
component $a$ remains constant.

Similarly, we define the rarefaction curve $\Rcal_2(U_0)$ issuing from
a given left-hand state $U_0$ and consisting of all the right-hand
states $U$ that can be connected to $U_0$ by a rarefaction wave
associated with the second characteristic field as
\begin{equation}
\Rcal_2(U_0):\quad
u=v_2(h,U_0):=u_0+2\sqrt{g}(\sqrt{h}-\sqrt{h_0}),\quad h\ge
h_0.\label{3.18}
\end{equation}
The $u$-component of the $2$-rarefaction wave is 
determined by \rf{3.17} and the $h$-component is given by
\rf{3.18}.

We can summarize the above results in: 

\begin{lemma}[Rarefaction wave curves]
Given a left-hand state $U_0$, the $1$-rarefaction curve
$\Rcal_1(U_0)$ consisting of all right-hand states $U$ that can be
connected to $U_0$ by a rarefaction wave associated with the first
characteristic field is 
$$
\Rcal_1(U_0):\quad
u=v_1(h,U_0):=u_0-2\sqrt{g}(\sqrt{h}-\sqrt{h_0}),\quad h\le h_0.
$$
The $2$-rarefaction curve $\Rcal_2(U_0)$ consisting of all
right-hand states $U$ that can be connected to $U_0$ by a
rarefaction wave associated with the second characteristic field is 
$$
\Rcal_2(U_0):\quad
u=v_2(h,U_0):=u_0+2\sqrt{g}(\sqrt{h}-\sqrt{h_0}),\quad h\ge h_0.
$$
\end{lemma}

We will also need backward curves which we define here for completeness.
Given a {\sl right-hand} state $U_0$, the $1$-rarefaction
curve $\Rcal_1^B(U_0)$ consisting of all {\sl left-hand} states $U$ that
can be connected to $U_0$ by a rarefaction wave associated with the
first characteristic field is 
\begin{equation}
\Rcal_1^B(U_0):\quad
u=v_1(h,U_0):=u_0-2\sqrt{g}(\sqrt{h}-\sqrt{h_0}),\quad h\ge h_0.
\label{3.19}
\end{equation}
The $2$-rarefaction curve $\Rcal_2^B(U_0)$ consisting of all
{\sl left-hand} states $U$ that can be connected to $U_0$ by a
rarefaction wave associated with the second characteristic field is 
\begin{equation}
\Rcal_2^B(U_0):\quad
u=v_2(h,U_0):=u_0+2\sqrt{g}(\sqrt{h}-\sqrt{h_0}),\quad h\le h_0.
\label{3.20}
\end{equation}

In turn, we are in a position to define the wave curves, as follows 
\begin{equation}
\aligned
&\WW_1(U_0)=\Scal_1(U_0)\cup\Rcal_1(U_0),\\
&\WW_1^B(U_0)=\Scal_1^B(U_0)\cup\Rcal_1^B(U_0),\\
&\WW_2(U_0)=\Scal_2(U_0)\cup\Rcal_2(U_0),\\
&\WW_2^B(U_0)=\Scal_2^B(U_0)\cup\Rcal_2^B(U_0).\\
\endaligned
 \label{3.21}
\end{equation}
Some properties of the wave curves are now checked. 

\begin{lemma}[Monotonicity properties] 
The wave curve $\WW_1(U_0)$ can be parameterized in the form $h\mapsto
u=u(h), {h>0}$, where the function $u$ is strictly convex and strictly decreasing
in $h$. The wave curve $\WW_2(U_0)$  can be parameterized in the form 
$h\mapsto u=u(h), h>0,$ where the function $u$ is strictly concave and strictly
decreasing in $h$.
\end{lemma}

\begin{proof} We only give the proof for the $1$-wave curve
$\WW_1(U_0)$, the proof for $\WW_2(U_0)$ being similar.
For the shock part $\Scal_1(U_0)$, we have
$$
{du\over dh} = -\sqrt{g\over 2}\dfrac{\dfrac{1}{
2h}+\dfrac{1}{h_0}+\dfrac{h_0}{2h^2}}{\sqrt{\dfrac{1}{h}+\dfrac{1}{h_0}}}<0.
$$
For the rarefaction part $\Rcal_1(U_0)$, we have
$$
{du\over dh}=-\sqrt{g\over h}<0.
$$
This establishes the desired monotonicity property of $\WW_1(U_0)$.

The convexity of $\WW_1(U)$ follows from the fact that $du/dh$ is
increasing. Indeed, along the shock part $\Scal_1(U_0)$ it holds
$$
{d^2u\over dh^2}=\sqrt{g\over 2}\dfrac{\Big({1\over
2h^2}+{h_0\over h^3}\Big)\sqrt{{1\over h}+{1\over h_0}} +{1\over
2h^2\sqrt{{1\over h}+{1\over h_0}}}\Big({1\over 2h}+{1\over
h_0}+{h_0\over 2h^2}\Big)}{{1\over h}+{1\over h_0}}>0
$$
and, along the rarefaction part $\Rcal_1(U_0)$, 
$$
{d^2u\over dh^2}= {\sqrt{g}\over 2h^{3/2}}>0, 
$$
which completes the proof.
\end{proof}

Next, we consider the $3$-curve from a state $U_0$, which consists of all states $U$ that can be connected to 
$U_0$ by a {\sl stationary wave.} As seen in \rf{1.9}, $U$ and $U_0$ are related
by the Rankine-Hugoniot relations
\begin{equation}
\aligned &[hu]=0\\
&[{u^2\over 2}+g(h+a)]=0.\\
\endaligned
 \label{3.22}
\end{equation}
This leads to a natural definition of a curve parameterized in $h$:
\begin{equation}
\WW_3(U_0):
\quad 
\begin{cases}
& u=u(h)={h_0u_0\over h},
\\
& a=a(h)=a_0+{u^2-u_0^2\over 2g}+h-h_0.
\end{cases} 
\label{3.23}
\end{equation}

%========================================================================================================

\section{Admissibility conditions for stationary waves}

\subsection{Two possible stationary jumps}

In view of the discussion in the previous section, the states across a stationary
wave are constraint by  the  Rankine-Hugoniot relations
\rf{3.22}. From a given left-hand state we have to determine the
right-hand state, which has three components, determined by the two
equations \rf{3.22}. Moreover, since the component $a$ changes
only through stationary waves (which propagate with zero speed) 
for given bottom levels $a_\pm$ we should solve for $u$ and $h$
in terms of $a$.
Thus, we rewrite \rf{3.22} in the form
$$
\aligned
 &u={h_0u_0\over h},\\
& a_0-a+{u^2-u_0^2\over 2g}+h-h_0=0. \endaligned
$$
Substituting for $u$ and re-arranging the terms, we obtain
$$
\aligned
 &u={h_0u_0\over h},\\
& a_0-a+{u_0^2\over 2g}\Big({h_0^2\over h^2}-1\Big)+h-h_0=0.
\endaligned
$$
This leads us to search for roots of the function
\begin{equation}
 \varphi(h) := a_0-a+{u_0^2\over 2g}\Big({h_0^2\over
 h^2}-1\Big)+h-h_0.
 \label{4.3}
\end{equation}
Let us set
\begin{equation}\aligned
& h_{\min}(U_0) :=\Big({u_0^2h_0^2\over g}\Big)^{1/3},\\
&a_{\min}(U_0) :=a_0+{u_0^2\over 2g}\Big({h_0^2\over
 h_{\min}^2}-1\Big)+h_{\min}-h_0.
 \endaligned
 \label{4.4}
\end{equation}

Some useful properties of the function $\varphi$ in \rf{4.3} are now derived. 

\begin{lemma} Suppose that $U_0 = (h_0, u_0, a_0)$ and $a$ are given with $u_0\ne 0$. 
The function $\varphi: (0, +\infty) \to \RR$ is smooth and convex and, for some $h_{\min}$,
it is decreasing in the interval $(0,h_{\min})$ and is increasing
in the interval $(h_{\min},\infty)$, with 
\begin{equation}
\lim_{h\to 0} \varphi(h)=\lim_{h\to +\infty} \varphi(h) = +\infty.
 \label{4.5}
\end{equation}
Furthermore, if $a\ge a_{\min}$ then the function $\varphi$ has two roots
$h_*(U_0), h^*(U_0)$ with $h_*(U_0)\le h_{\min}(U_0)\le
h^*(U_0)$. These inequalities are strict whenever $a>a_{\min}(U_0)$.
\end{lemma}

\begin{proof} The smoothness of the function $\varphi$  and the limiting conditions are obvious.
Moreover, we have 
$$
{d\varphi(h)\over dh} = -{u_0^2h_0^2\over gh^3}+1
$$
(for $u_0\ne 0$) which is positive if and only if
$$
h > \Big({u_0^2h_0^2\over g}\Big)^{1/3}=h_{\min}(U_0).
$$
This establishes the monotonicity property of $\varphi$.
Furthermore, we have 
$$
{d^2\varphi(h)\over dh^2}={3u_0^2h_0^2\over gh^4}\ge 0,
$$
which shows the convexity of $\varphi$. If
$a>a_{\min}(U_0)$, then $\varphi(h_{\min}(U_0))<0$.  The other
conclusions follow immediately.
\end{proof}

It is straightforward to check: 

\begin{lemma} The function $h_{\min}$ satisfies the following inequalities: 
\begin{equation}
\aligned
& h_{\min}(U_0) >h_0,\quad\quad U_0\in A_1\cup A_3,\\
& h_{\min}(U_0) <h_0,\quad\quad U_0\in A_2,\\
& h_{\min}(U_0) =h_0,\quad\quad U_0\in \CC,\\
 \endaligned
 \label{4.6}
\end{equation}
The roots $h^*$ and $h_*$ satisfy the following inequalities: 
 \begin{itemize}
 \item[(i)]
 If $a > a_0$, then
\begin{equation}
h_*(U_0) < h_0 < h^*(U_0).
 \label{4.7}
\end{equation}
\item[(ii)] If $a < a_0$, then
\begin{equation}
\aligned
&h_0 < h_*(U_0)  \quad\quad U_0 \in A_1\cup A_3,\\
&h_0 > h^*(U_0)  \quad\quad U_0 \in A_2.
\endaligned
 \label{4.8}
\end{equation}
 \end{itemize}
The function $a_{\min}(U_0)$ satisfy the following inequalitites:  
\begin{equation}
\aligned
&a_{\min}(U_0) < a_0,\quad (h_0,u_0)\in A_1 \cup A_2 \cup A_3, 
\\
&a_{\min}(U_0) = a_0,\quad (h_0,u_0)\in \CC_\pm.  
\endaligned
\label{4.9}
\end{equation}
\end{lemma}

The states that can be connected by stationary waves are characterized as follows. 

\begin{proposition} 
\label{theo41} 
Fix a left-hand state $U_0=(h_0,u_0,a_0)$ and a right-hand bottom level $a$.
\begin{itemize}
\item[(i)] If $u_0\ne 0$ and $a>a_{\min}(U_0)$, then there are two
distinct right-hand states
$$
U_{1,2}:=(h_{1,2}(U_0),u_{1,2}(U_0),a)
$$
where $u_i(U_0):=h_0u_0/h_i(U_0), i=1,2$, that can be connected to
$U_0$ by a stationary wave satisfying the Rankine-Hugoniot
relations.
\item[(ii)] If $u_0\ne 0$ and $a=a_{\min}(U_0)$, the two states in
(i) coincide and we obtain a unique stationary wave.
\item[(iii)] If $u_0\ne 0$ and $a<a_{\min}(U_0)$, there is no
stationary wave from $U_-$ to a state with level $a$.
\item[(iv)] If $u_0=0$, there is only one stationary jump defined
by
$$
u=u_0=0,\quad h=h_0+a-a_0.
$$
\end{itemize}
\end{proposition}

We arrive at an important conclusion on stationary jumps.

\begin{proposition} \label{prop42} 
For $u_0\ne 0$, the state $(h_{1}(U_0),u_{1}(U_0))$ belongs to $A_1$ if $u_0 < 0$,
and belongs to $A_3$ if $u_0 > 0$, while the state
$(h_{2}(U_0),u_{2}(U_0))$ always belongs to $A_2$.
Moreover, we have
\begin{equation}
(h_{\min}(U_0), u=h_0u_0/h_{\min}(U_0)) \in 
\begin{cases}
\CC^+,         & u_0 > 0,
\\
\CC^-,         & u_0 < 0. 
\end{cases} 
\label{4.10}
\end{equation}
\end{proposition}

It is interesting to observe that the shock speed in the genuinely nonlinear
characteristic fields will change sign along the shock curves.
Therefore, it exchanges its order with the linearly degenerate
field, as stated in the following theorem.

\begin{proposition} 
\label{theo43}
  (a) If  $U_0\in A_1$, then there exists $\tilde U_0\in\Scal_1(U_0)\cap  A_{2}^+$
   corresponding to $h=\tilde h>h_0$ such that
\begin{equation}
\aligned
&\bar\ld_1(U_0,\tilde U_0) = 0,\\
&\bar\ld_1(U_0,U) > 0,\quad  U\in\Scal_1(U_0), h \in (h_0,\tilde h_0),\\
&\bar\ld_1(U_0,U) < 0,\quad U\in\Scal_1(U_0), h \in (\tilde h_0,+\infty).\\
\endaligned
 \label{4.11}
\end{equation}
If $U_0\in A_2\cup A_3$, then
\begin{equation}
\bar\ld_1(U_0,U) < 0, \quad  U\in \Scal_1(U_0).
 \label{4.12}
\end{equation}

(b) If  $U_0\in A_3$, then there exists $\bar U_0\in\Scal_2^B(U_0)\cap A_{2}^-$ corresponding to $h=\bar h>h_0$ such that
\begin{equation}
\aligned
&\bar\ld_2(U_0,\bar U_0) = 0,\\
&\bar\ld_2(U_0,U) > 0,\quad U\in\Scal_2^B(U_0), h \in (h_0,\bar h_0),\\
&\bar\ld_2(U_0,U) < 0,\quad U\in\Scal_2^B(U_0), h \in (\bar h_0,+\infty).\\
\endaligned
 \label{4.13}
\end{equation}

If $U_0\in A_1\cup A_2$, then
\begin{equation}
\bar\ld_2(U_0,U) > 0,\qquad  U\in \Scal_2^B(U_0).
 \label{4.14}
\end{equation}
\end{proposition}

%-----------------------------------------------------------------------

\subsection{Two-parameter wave sets}

From Proposition~\ref{prop42} and the arguments in the previous section, we can now construct
wave composites. It turns out that two-parameter wave sets can be constructed. 
For definiteness, we now illustrate this feature on a particular case.
Suppose that $U_0=(h_0,u_0,a_0)\in A_2^+$. We can use a stationary wave from $U_0$
to a state $U_m=(h_m,u_m,a_m)\in A_2^+$ using $h^*$, followed by
another stationary wave from $U_m$ to $U\in A_1$ using the
corresponding value $h_*$, then we continue with $1$-waves, and as in $A_1$ the
characteristic speed is positive. As $a_m$ can vary, the set of such states 
$U$ forms a two-parameter set of composite waves containing first and third waves.
Such wave sets were constructed even for strictly hyperbolic systems by Hayes and LeFloch \cite{HayesLeFloch}. 

To make the Riemann problem well-posed, it is necessary to impose an additional admissibility criterion.

\subsection{The monotonicity criterion}

Since the Riemann problem for \rf{1.1} may in principle 
admit up to a one-parameter family of solutions, 
we now require that the Riemann solutions of interest satisfy a monotonicity condition in the component $a$.

\begin{itemize}

 \item[(MC)] (Monotonicity Criterion) \quad Along any
stationary curve $\WW_3(U_0)$, the bottom level $a$ is 
a monotone function in $h$.  The total variation of the bottom level
component of any Riemann solution must not exceed (and, therefore,
is equal to) $|a_L-a_R|$, where $a_L, a_R$ are left-hand and
right-hand cross-section levels.
\end{itemize}

A similar selection criterion was used by Isaacson and Temple
\cite{IsaacsonTemple92, IsaacsonTemple95} and by LeFloch and Thanh
\cite{LeFlochThanh03.2}, and by Goatin and LeFloch \cite{GoatinLeFloch}.
Under the transformation (if necessary) 
$$
x \to -x, \qquad u \to -u,
$$
a right-hand state $U=(h,u,a)$ transforms into a left-hand state of
the form $U'=(h,-u,a)$. Therefore, it is not restrictive to assume that
\begin{equation}
a_L < a_R.
 \label{5.1}
\end{equation}

\begin{lemma} \label{lem51} The Monotonicity Criterion implies that stationary shocks
do not cross the boundary of strict hyperbolicity. In other
words, we have: 
\begin{itemize}
 \item[(i)] If $U_0\in A_1\cup A_3$, then
only the stationary shock based on the value $h_*(U_0)$ is admissible.
\item[(ii)] If $U_0\in A_2$, then only the stationary shock using $h^*(U_0)$ is admissible.
\end{itemize}
\end{lemma}

%-----------------------------------------------------------------------------

\begin{proof}
Recall that the Rankine-Hugoniot relations associated with the linearly
degenerate field \rf{3.23} implies that the component $a$ can be
expressed as a function of $h$:
$$
 a=a(h)=a_0+{u^2-u_0^2\over 2g}+h-h_0,
$$
where
$$
u=u(h)={h_0u_0\over h}.
$$
Thus, differentiating $a$ with respect to $h$, we find
$$
\aligned
a'(h)&={uu'(h)\over g}+1=-u{h_0u_0\over gh^2}+1\\
&=-{u^2\over  gh}+1\\
\endaligned
$$
which is positive (resp. negative) if and only if
$$
u^2 - g \, h < 0 \qquad \text{(resp. $u^2 - g \, h > 0$)} 
$$
or $(h,u,a)\in A_2$ (resp. $\in A_1$ or $\in A_3$). Thus, in order that
$a'$ keeps the same sign, the point $(h,u,a)$ must remain on the
same side as $(h_0,u_0,a_0)$ with respect to $\CC_\pm$. The
conclusions in (i) and (ii) follow.
\end{proof}

It follows from Lemma \ref{lem51} that for a given $U_0=(h_0,u_0,a_0)\in A_i, i=1,2,3,$ and a level $a$, we can
define a unique point $U=(h,u,a)$ so that the two points $U_0,U$
can be connected by a stationary wave satisfying the (MC)
criterion. We have a mapping 
\begin{equation}
\aligned
 SW(.,a): [0,\infty)\times \RR\times \RR_+ & \to [0,\infty)\times \RR\times\RR_+\\
U_0=(h_0,u_0,a_0) &\mapsto  SW(U_0,a)= U = (h,u,a),\\
\endaligned
 \label{5.2}
\end{equation}
such that $U_0$ and $U$ can be connected by a stationary wave
satisfying the (MC) condition. Observe that this mapping is
single-valued except on the hypersurface $\CC$, where it has
two-values.

Let us use the following notation: $W_i(U_0,U)$ will stand for the
$i$th wave from a left-hand state $U_0$ to the right-hand state
$U$, $i=1,2,3$. To represent the fact that the wave $W_i(U_1,U_2)$
is followed by the wave $W_j(U_2,U_3)$, we use the notation:
\begin{equation}
W_i(U_1,U_2) \oplus W_j(U_2,U_3).
 \label{5.3}
\end{equation}

%==========================================================================================================

\section{The Riemann problem}

In this section we construct the solutions of the Riemann problem, by combining 
Lax shocks, rarefaction waves, and stationary waves
satisfying the admissibility condition (MC).

Recall that for general strictly hyperbolic systems of conservation laws, 
the solution to the Riemann problem exist when the initial jump is sufficiently small only. 
That is to say that the right-hand states $U_R$ should lie in a small neighborhood of the
left-hand state $U_L$.  However, for the system \rf{1.1}, we can cover large data 
and essentially cover a full domain of existence for any given
left-hand state. More precisely, we determine the precise range of
right-hand states in which the Riemann solution exists.

\subsection{Solutions containing only one wave of each characteristic family}   

We begin by constructing solutions containing only one
wave corresponding to each characteristic field which is
identified as each family of waves. This structure of solutions is standard
 in the theory of {\sl strictly hyperbolic} system of conservation
laws. In the next subsection we will consider solutions that contain
up to two waves in the same family. The following theorem deals
with the case where the left-hand state $U_L$ is in $A_1$.

\begin{figure}   \includegraphics[width=0.7\textwidth]{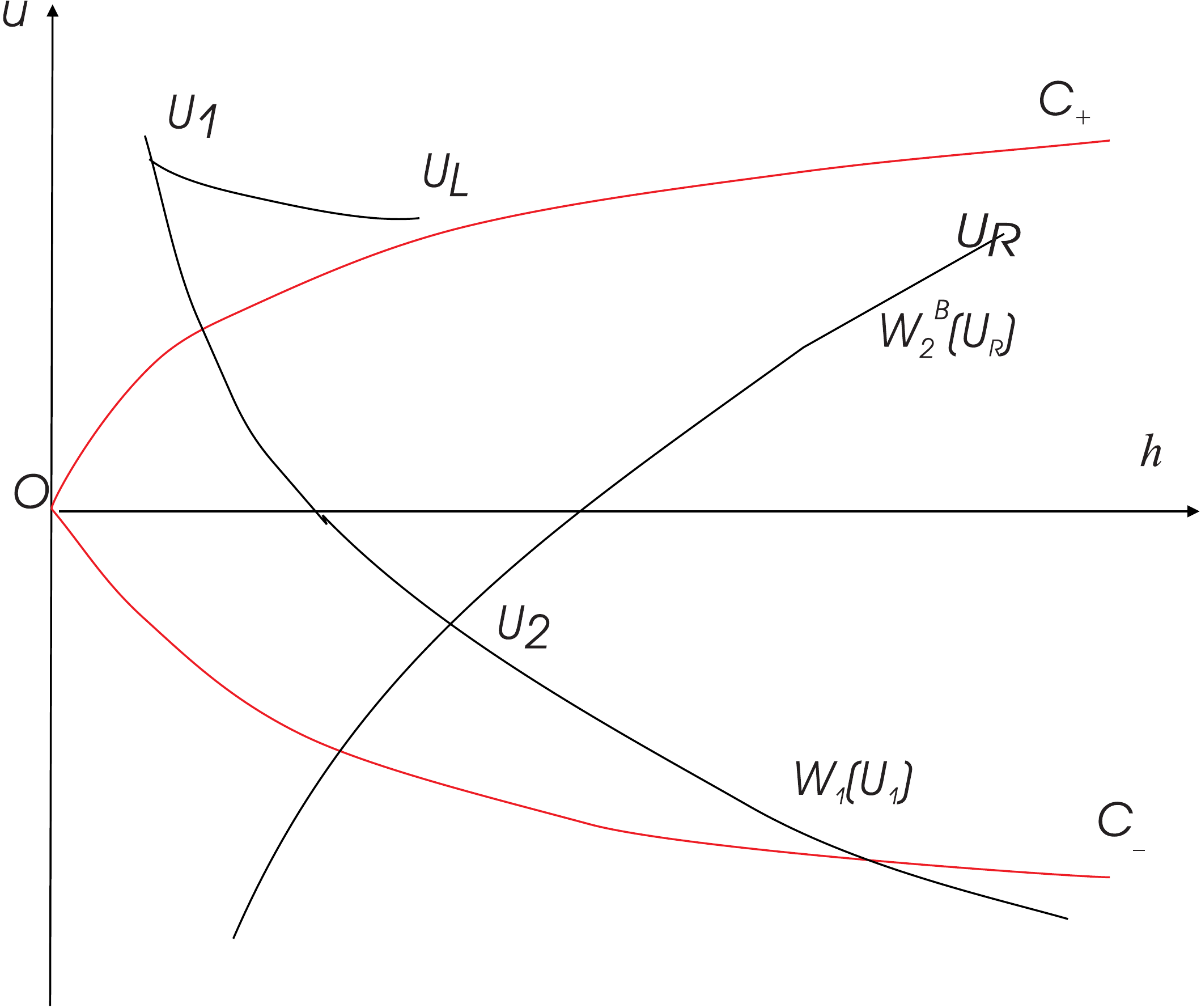}\\
  \caption{Solution for $U_L\in A_1$}\label{fig41}
\end{figure}

\begin{theorem} \label{theo51}
 Let $U_L\in A_1$ and set $U_1:=SW(U_L,a_R)$,
$\{U_2\}=\WW_1(U_1)\cap \WW_2^B(U_R)$. Then, the
 Riemann problem \rf{1.1}-\rf{1.2} admits an admissible solution with the following structure
\begin{equation}
W_3(U_L,U_1) \oplus W_1(U_1,U_2) \oplus W_2(U_2,U_R),
 \label{5.4}
\end{equation}
provided $h_2\le \tilde h_1$. (Figure \ref{fig41}).
\end{theorem}

\begin{proof} Observe that the set of composite waves $SW(\WW_1(U_L), a_R)$
consists of three monotone decreasing curves, and each lies entirely
in each region $A_i, i=1,2,3$.  The monotone increasing backward
curve $\WW_2^B(U_R)$ therefore may cut the three composite curves
at a unique point, two point, or else does not meet the wave composite
set. The Riemann problem therefore may admit a unique solution,
two solutions, or has no solution.

The state $U_L$ belongs to $A_1$ and in this region, the $\ld_3$
is the smallest of the three characteristic speeds. A stationary
wave from $U_L=(h_L,u_L,a_L)$ to $U_1=(h_1,u_1,a_R)$ exists, since
$a_L\le a_R$. Moreover, by Lemma~\ref{lem51}, we have $U_1\in A_1$.

If $h_2\le h_1$, then the stationary wave is followed by a
$1$-rarefaction wave with positive speed, and then can be
continued by a $2$-wave $W_2(U_2,U_R)$. If $h_2>h_1$, then the
$1$-wave in \rf{5.3} is a shock wave. Since $h_2\le \bar h_1$ and
$U_1\in A_1$, the shock speed $\ld_2(U_1,U_2)\ge 0$, and thus it
can follow a stationary wave (with zero speed). Moreover, it is
derived from \rf{3.5} that
$$
\aligned
\bar\ld_{1}(U_1,U_2) &=u_1-\sqrt{{g\over 2}\Big(h_2+{h_2^2\over h_1}\Big)},\\
&={h_2u_2-h_1u_1\over h_2-h_1}\\
&=u_2-\sqrt{{g\over 2}\Big(h_1+{h_1^2\over h_2}\Big)},
\endaligned 
$$
thus 
\begin{equation}
\aligned
\bar\ld_{1}(U_1,U_2)
&\le u_2+\sqrt{{g\over 2}\Big(h_R+{h_R^2\over h_2}\Big)}
=\bar\ld_{2}(U_2,U_R).
\endaligned \label{5.5}
\end{equation}
This means the $1$-shock $S_1(U_1,U_2)$ can always follow the
$2$-shock $ S_2(U_2,U_R)$. Similar for rarefaction waves.
Therefore, the solution structure \rf{5.3} holds.
\end{proof}

The following theorem deals with the case where the left-hand state $U_L$ is in $A_1\cup A_2$.

\begin{theorem}
\label{theo52}
Let $U_L\in A_1\cup A_2$. Then there exists a region of values $U_R$
such that $SW(\WW_1(U_L),a_R)\cap\WW_2^B(U_R)\ne \emptyset$. In
this case this intersection may contain  either only one or both
points $U_1\in A_2$ and $U_2\in A_3$. The Riemann problem
\rf{1.1}-{1.2} therefore has a solution with the structure
\begin{equation}
W_1(U_L,U_3) \oplus W_3(U_3,U_1) \oplus W_2(U_1,U_R),
 \label{5.6}
\end{equation}
where $U_3\in \WW_1(U_L)$ is the point such that
$U_1=SW(U_3,a_R)$, and also 
\begin{equation}
W_1(U_L,U_4) \oplus W_3(U_4,U_2) \oplus W_2(U_2,U_R),
 \label{5.7} 
\end{equation}
where $U_4\in  \WW_1(U_L)$ is the point such that
$U_2=SW(U_4,a_R)$, if $h_2\ge \bar h_R$ whenever $U_R\in A_2^-$.
(Figure 3 ) 
\end{theorem}

\begin{figure}
  \includegraphics[width=0.7\textwidth]{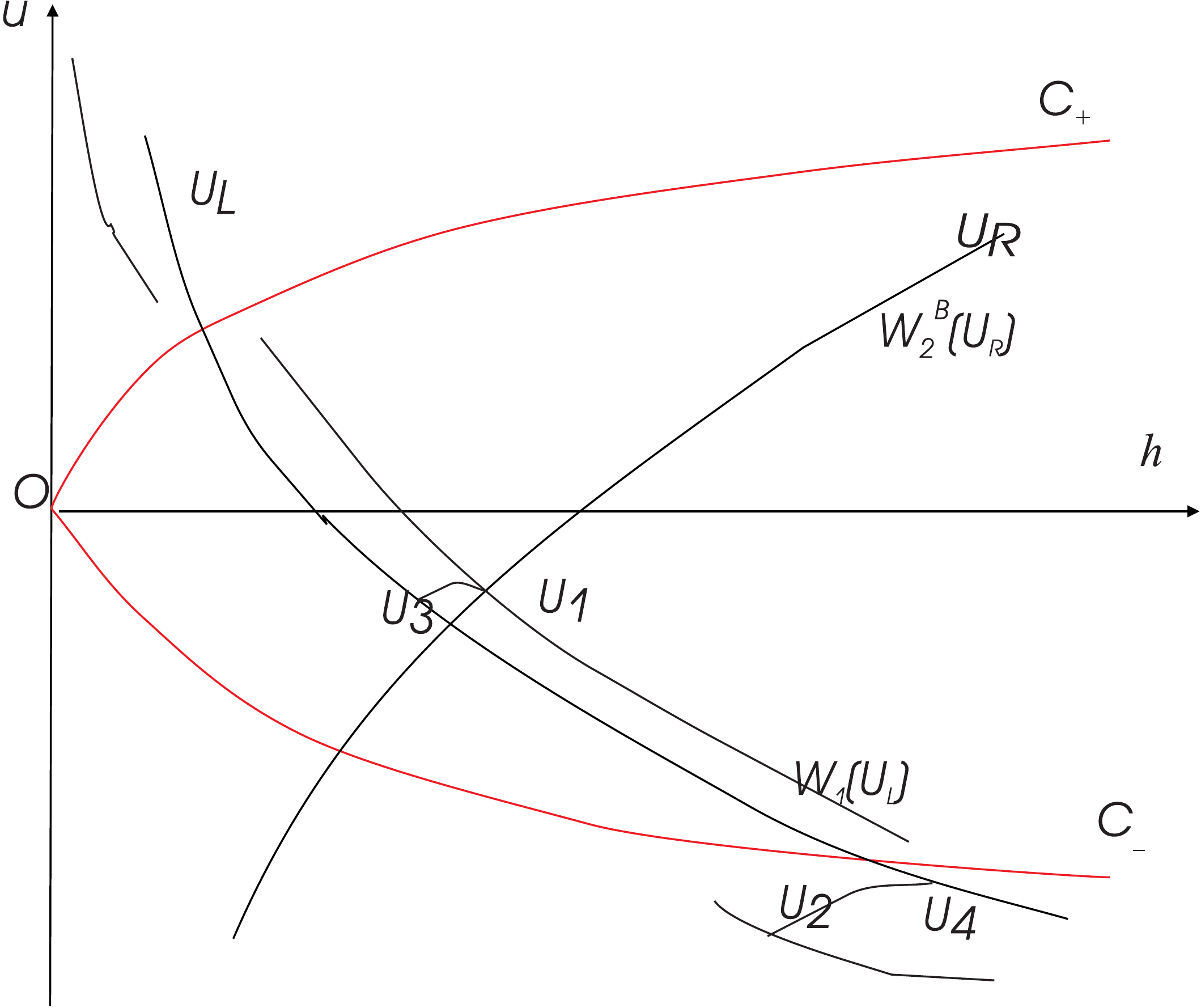}\\
  \caption{Solution for $U_L\in A_1\cup A_2$}\label{fig42}
\end{figure}

\begin{proof}
 The solution may begin with a $1$-wave, either $1$-shock with a negative shock
speed to a state $U_3$, or a $1$-rarefaction wave with
$\ld_1(U_3)\le 0$, followed by a stationary wave $W_3(U_3,U_1)$
from $U_3$ to $U_1$, then followed by a $2$-wave $W_2(U_1,U_R)$
from $U_1$ to $U_R$. It is similar in the case of $U_2$. However,
in order that the stationary wave $W_3(U_4,U_2)$, for some $U_4\in
\WW_1(U_L)$ and $U_4\in A_3$ obviously, to be followed by a
$2$-wave $W_2(U_2,U_R)$, it is required that the wave is a shock
with non-negative shock speed $\ld_2(U_2,U_R)$. This is
equivalent to $h_2\ge \bar h_R$.
\end{proof}

\begin{figure}
  \includegraphics[width=0.7\textwidth]{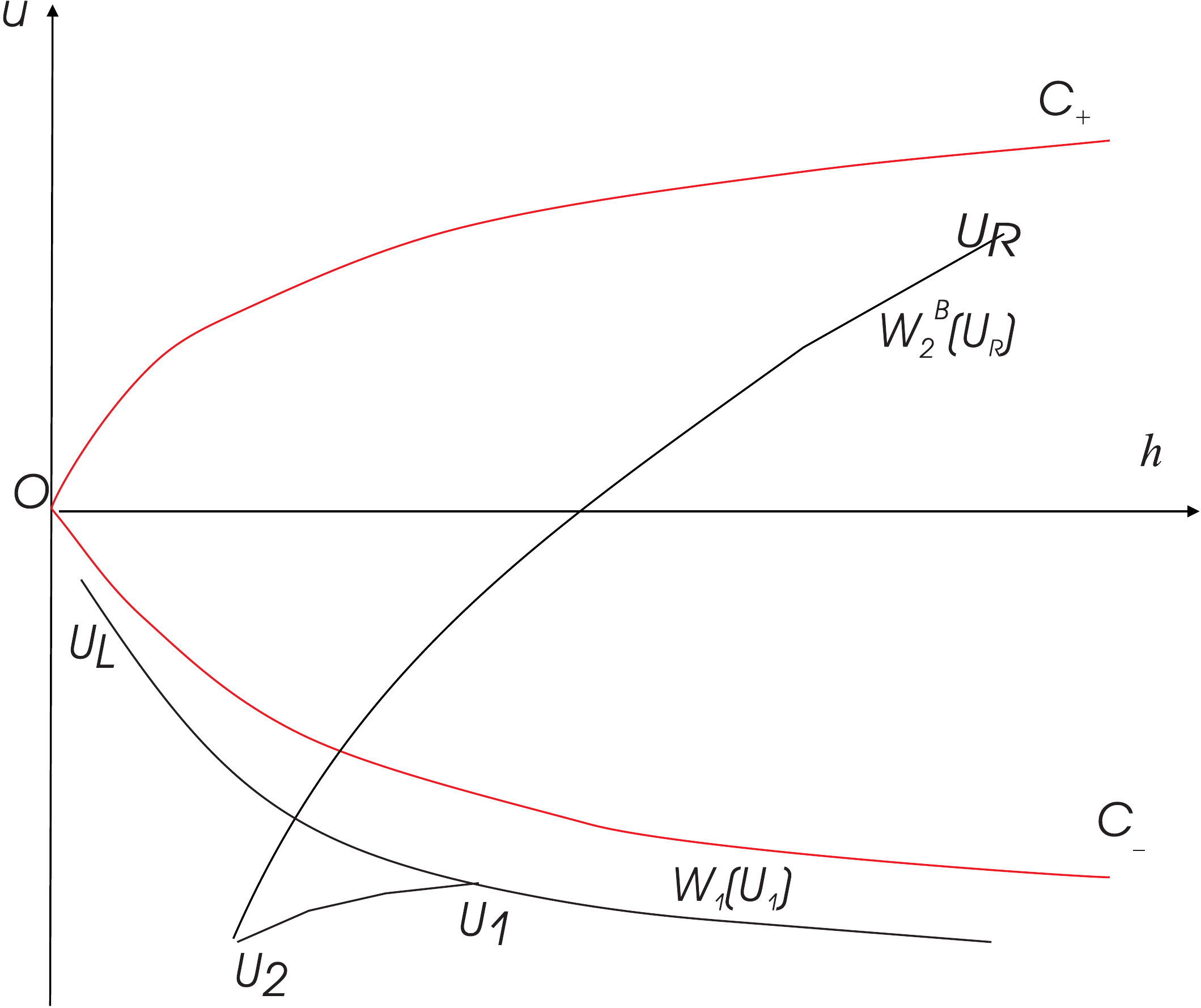}\\
  \caption{Solution for $U_L\in A_3$}\label{fig43}
\end{figure}

\begin{theorem} \label{theo53} Let $U_L\in A_3$ and $U_R\in A_1\cup A_2$, and 
set $U_1=SW(\WW_2^B(U_R),a_L)\cap \WW_1(U_L)$ and 
$U_2=SW(U_1,a_R)\in \WW_2^B(U_R)$.
\begin{itemize}
\item[(i)] If $U_1\in A_2^+\cup \CC_+\cup\{u=0\}$, the Riemann
problem \rf{1.1}-\rf{1.2} has a solution with the following structure
\begin{equation}
W_1(U_L,U_1) \oplus W_3(U_1,U_2) \oplus W_2(U_2,U_R).
 \label{5.8}
\end{equation}

\item[(ii)] If $U_1\in A_2^-\cup \CC_-$, provided $h_R\ge \bar
h_2$, the Riemann solution \rf{5.8} also exists.

 \item[(iii)] If $U_1\in A_1\cup A_3$, the
construction \rf{5.8} does not make sense.
\end{itemize}
(Figure \ref{fig43}).
\end{theorem}

\begin{proof} If $U_1\in A_2\cup\CC$, the non-positive speed wave
$W_1(U_L,U_1)$ can be followed by a stationary wave
$W_3(U_1,U_2)$. 

 When $U_2\in A_2-$, if $U_1\in A_2^+\cup \CC_+\cup\{u=0\}$, then
 this stationary wave can always be followed by a $2$-wave
 $W_2(U_2,U_R)$, since the wave speed of the $2$-wave is positive.
 This establishes (i).

 If   $U_1\in A_2^-\cup \CC_-$, the wave speed of the $2$-wave $W_2(U_2,U_R)$ is non-negative
 if and only if $h_R\ge \bar h_2$. This proves (ii).

 If $U_1\in A_1$, the $1$-wave has positive speed. So it can not
 be followed by a stationary wave. If $U_1\in A_3$, then $U_2\in A_3$ by the (MC) criterion.
 So the  $2$-wave $W_2(U_2,U_R)$ has negative speed. So it can not
 be proceeded by a stationary wave. This proves (iii).
\end{proof}

The above theorem enables $U_R$ to vary in each region $A_1, A_2$, 
and $A_3$. The next theorem enables $U_L$ to vary in all the three regions.

\begin{figure}
  \includegraphics[width=0.7\textwidth]{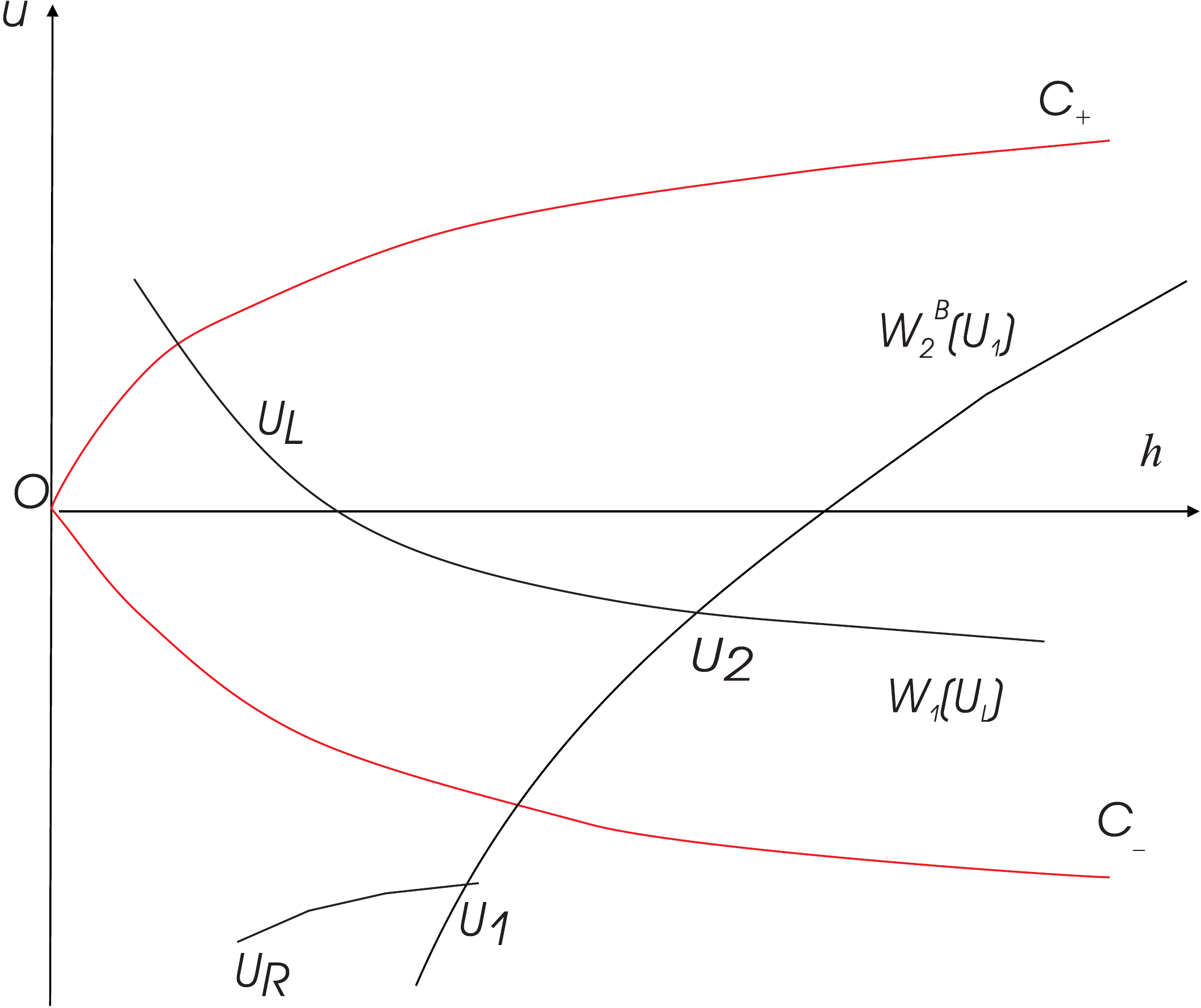}\\
  \caption{$U_L$ may be anywhere}\label{fig44}
\end{figure}

\begin{theorem} \label{theo54}
 Let $U_R\in A_3$. Set $U_1=SW(U_R,a_L)$,
$U_2=\WW_2^B(U_1)\cap \WW_1(U_L)$. A Riemann solution exists and
has the following structure

\begin{equation}
W_1(U_L,U_2) \oplus W_2(U_2,U_1) \oplus W_3(U_1,U_R),
 \label{5.9}
\end{equation}
provided $ h_2\le \bar h_1$. (Figure \ref{fig44}).
\end{theorem}

\begin{proof} The stationary wave $W_3(U_1,U_R)$ turns out to have
the greatest wave speed. In order for this wave to be proceeded by
the $2$-wave $W_2(U_1,U_2)$, the wave speed of this $2$-wave has
to be non-positive. This is equivalent to the  condition $h_2\le
\bar h_1$, according to Proposition~\rf{theo43}. Similar to \rf{5.7},
we have
$$
\ld_1(U_L,U_2)\le \ld_2(U_2,U_1).
$$
 so that the $1$-wave $W_1(U_L,U_2)$ can follow the
$2$-wave $ W_2(U_2,U_1)$.
\end{proof}

%---------------------------------------------------------------------------------------------------------

\subsection{Solutions containing more than one wave of each characteristic family}

It is remarkable feature of the shallow water system that we can
also construct solutions with {\sl four} elementary waves, 
using three available characteristic fields. This illustrates one of the difficulties  
in coping with the Riemann problem when the system under consideration is not strictly hyperbolic.

\begin{figure}
  \includegraphics[width=0.7\textwidth]{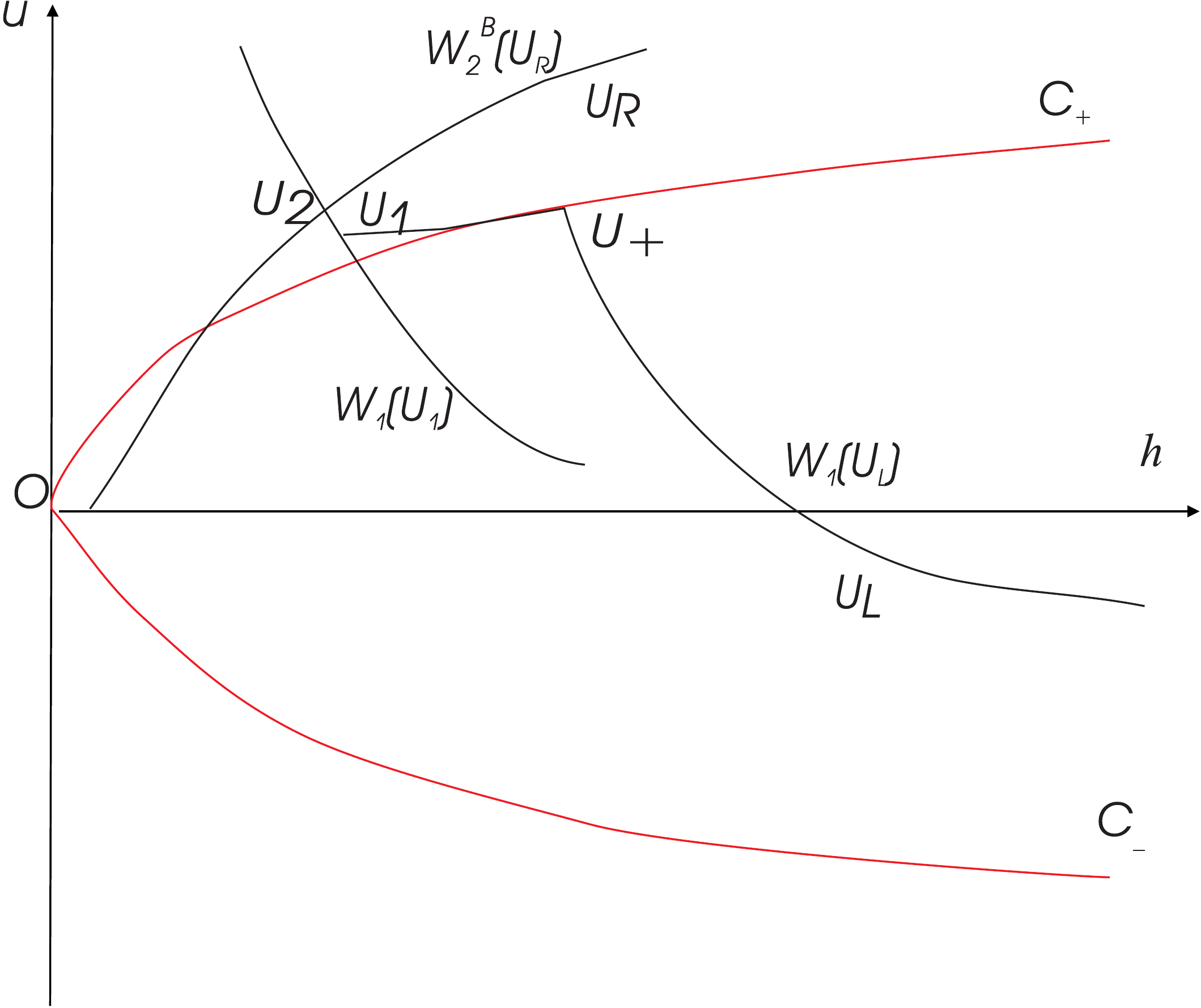}\\
  \caption{Solution with repeated two $1$-waves}\label{fig45}
\end{figure}

\begin{theorem} Let $U_L\in A_2\cup A_3$ and set
$U_+=\WW_1(U_L)\cap \CC_+,\{U_1\}=SW(U_+,a_R)\cap
A_1,\{U_2\}=\WW_1(U_1)\cap\WW_2^B(U_R)$.  
The Riemann problem
\rf{1.1}-\rf{1.2} has a solution with the following structure
\begin{equation}
R_1(U_L,U_+) \oplus W_3(U_+,U_1) \oplus W_1(U_1,U_2)\oplus
W_2(U_2,U_R),
 \label{5.11}
\end{equation}
provided $h_2\le \tilde h_1$. (Figure \ref{fig45}).
\end{theorem}

\begin{figure}
  \includegraphics[width=0.7\textwidth]{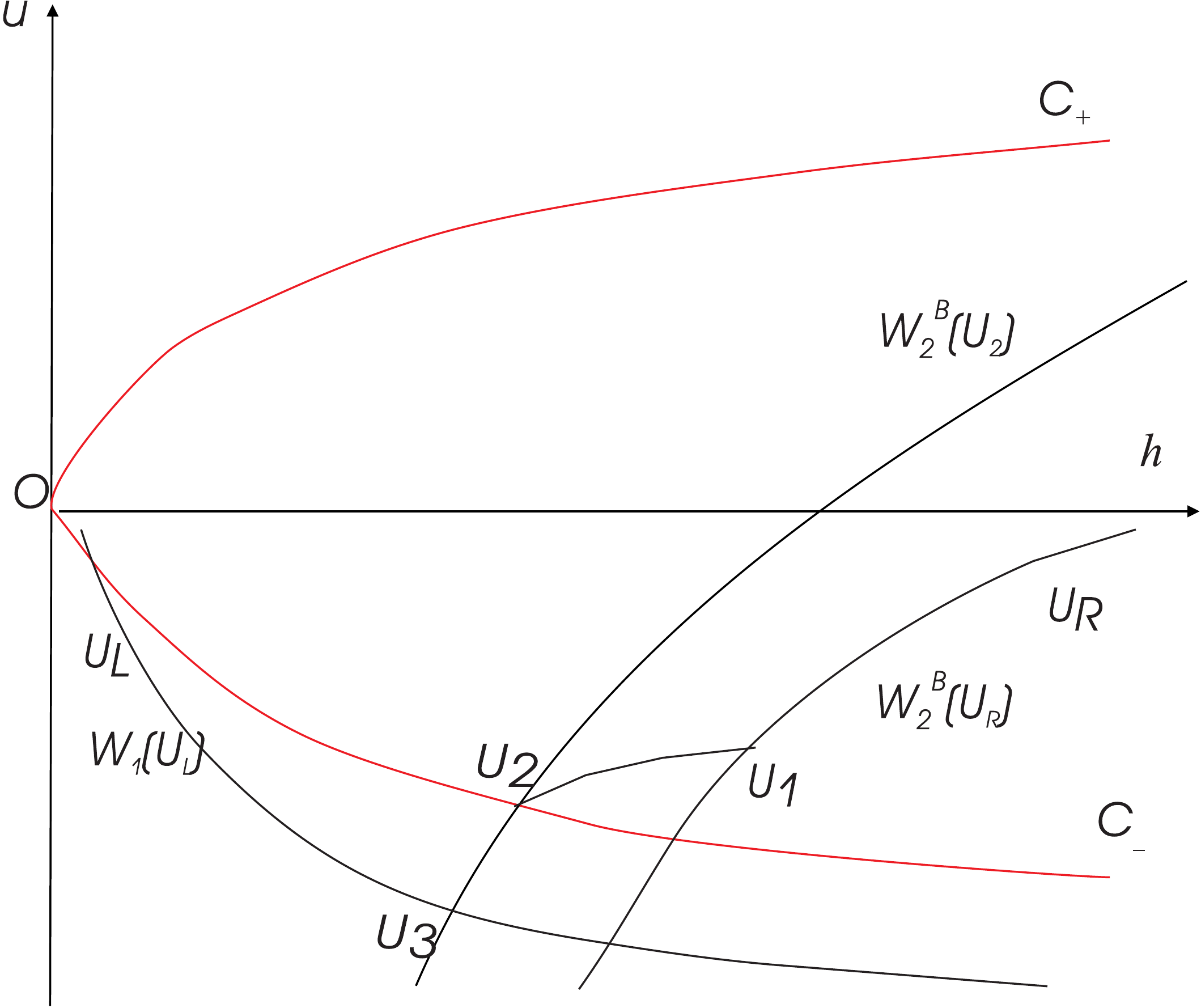}\\
  \caption{Solution with repeated two $2$-waves}\label{fig46}
\end{figure}

\begin{theorem} For any $U_L$, set $ \{U_1\}=SW(\CC_-,a_R)\cap  \WW_2^B(U_R)\cap A_2$, 
$U_2=(h_2,u_2,a_L)\in \CC_-$ such that $U_1=SW(U_2)$, and 
$\{U_3\}=\WW_2^B(U_2)\cap \WW_1(U_L).$ Then the Riemann problem
\rf{1.1}-\rf{1.2} has a solution with the following structure
\begin{equation}
W_1(U_L,U_3) \oplus R_2(U_3,U_2) \oplus W_3(U_2,U_1)\oplus
W_2(U_1,U_R),
 \label{5.12}
\end{equation}
provided $h_R\ge \bar h_1$ and $h_3\le h_2$. (Figure \ref{fig46}).
\end{theorem}

Thus, we see that the Riemann problem
\rf{1.1}-\rf{1.2} has a solution consisting of a $1$-, a $3$-, and
two $2$-waves.

It is interesting to note that there are solutions satisfying the
(MC) criterion which contain three waves with the same speed (zero).
This is the case when a stationary wave jumps from the level
$a=a_L$ to an intermediate level $a_m$ between $a_L$ and $a_R$,
followed by an "intermediate" $k$-shock with zero speed at the
level $a_m$, $k=1,3$, and then followed by  another stationary
wave jumping from the level $a_m$ to $a_R$. Thus, there are only
two possibilities:
\begin{itemize} 
\item[(i)]
$U_L$ belongs to $A_1$ and a $1$-shock with zero speed is used.
\item[(ii)]
 $U_R$ belongs to $A_3$ and a $2$-shock with zero speed is used.
\end{itemize}

We just describe the first case (i), as the second case is similar.
Recall from Proposition~\ref{theo43} that for any $U\in A_1$, there
exists a unique point denoted by $\tilde U\in \WW_1(U)\cap A_2$
such that
$$
\bar\ld_1(U,\tilde U)=0.
$$

\begin{theorem} Let $U_L\in A_1$ and set
$$
\aligned
&SW(U_L,[a_L,a_R]) :=\cup_{a\in [a_L,a_R]} SW(U_L,a),\\
&\widetilde{SW}(U_L,[a_L,a_R]):=\{\tilde U | \ U \in
SW(U_L,[a_L,a_R])\}.\\
\endaligned
$$
Whenever
$$
\widetilde{SW}(U_L,[a_L,a_R])\cap \WW_2^B(U_R)\ne \emptyset
$$
there exist $a_m\in [a_L,a_R]$, $U_1=SW(U_L,a_m)$, and 
$$
U_2\in \widetilde{SW}(U_L,[a_L,a_R])\cap \WW_2^B(U_R)
$$
that defines a solution with the structure
\begin{equation}
W_3(U_L,U_1) \oplus S_1(U_1,\tilde U_1) \oplus W_3(\tilde
U_1,U_2)\oplus W_2(U_2,U_R).
 \label{5.14}
\end{equation}
\end{theorem}

%================================================================================================================= 

\section*{Acknowledgments}

The first author (P.G.L.) was supported by the A.N.R. Grant 06-2-134423: 
{\em Mathematical methods in general relativity} (MATH-GR) and the Centre National de la Recherche Scientifique (CNRS).

%====================================================================================================================

\end{document}